\documentclass[12pt]{article}
\def\Dj{\hbox{D\kern-.73em\raise.30ex\hbox{-}
\raise-.30ex\hbox{}}}
\def\dj{\hbox{d\kern-.33em\raise.80ex\hbox{-}
\raise-.80ex\hbox{\kern-.40em}}}
\usepackage{amsmath,amssymb,amsthm, amscd, amsfonts}
\allowdisplaybreaks[4]
\newtheorem {Lemma}{Lemma}

\newtheorem{Proposition}{Proposition}

\newenvironment {Proof} {\noindent {\bf Proof.}}{\hspace*{\fill}$\Box$\par\vspace{3mm}}

\begin{document}

\baselineskip=0.31in

\vspace*{40mm}

\begin{center}
{\Large \bf ON SUM OF POWERS OF LAPLACIAN EIGENVALUES AND LAPLACIAN
ESTRADA INDEX OF GRAPHS}

\vspace{10mm}

{\large \bf Bo Zhou}

\vspace{6mm}

\baselineskip=0.20in

{\it Department of Mathematics, South China Normal University,
\\
Guangzhou 510631, P. R. China\/}\\
{\rm e-mail:} {\tt zhoubo@scnu.edu.cn}

\vspace{6mm}

(Received June 30, 2008)
\end{center}

\vspace{6mm}

\baselineskip=0.27in

\noindent {\bf Abstract }

\vspace{4mm}

Let $G$ be a simple graph and $\alpha$ a real number. The quantity
$s_{\alpha}(G)$ defined as the sum of the $\alpha$-th power of the
non-zero Laplacian eigenvalues of $G$ generalizes several concepts
in the literature. The Laplacian Estrada index is a newly introduced
graph invariant based on Laplacian eigenvalues. We establish bounds
for $s_{\alpha}$ and Laplacian Estrada index related to the degree
sequences.

\baselineskip=0.25in

\vspace{6mm}

\begin{center}
{\bf 1. INTRODUCTION}
\end{center}

Let $G$ be a simple graph possessing $n$ vertices. The Laplacian
spectrum of $G$, consisting of the numbers $\mu_1, \mu_2, \ldots,
\mu_n$ (arranged in non-increasing order), is the spectrum of the
Laplacian matrix  of $G$.
It is known that
$\mu_n=0$ and the multiplicity of $0$ is equal to the number of
connected components of $G$. See \cite{Mer,Moh} for more details for
the properties of the Laplacian spectrum.

Let $\alpha$ be a real number and let $G$ be a graph with $n$
vertices. Let $s_{\alpha}(G)$ be the sum of the $\alpha$-th power of
the non-zero Laplacian eigenvalues of $G$, i.e.,
\[
s_{\alpha}(G)=\sum_{i=1}^h\mu_i^{\alpha},
\]
where $h$ is the number of non-zero Laplacian eigenvalues of $G$.
The cases  $\alpha=0, 1$ are trivial as $s_0(G)=h$ and $s_1(G)=2m$,
where $m$ is the number of edges of $G$. For a nonnegative integer
$k$, $t_k(G)=\sum_{i=1}^n\mu_i^{k}$ is the $k$-th Laplacian spectral
moment of $G$. Obviously, $t_0(G)=n$ and $t_k(G)=s_k(G)$ for $k\ge
1$.  Properties of $s_{2}$ and $s_{\frac{1}{2}}$ were studied
respectively in \cite{Laz} and \cite{LL}. For a connected graph $G$
with $n$ vertices, $ns_{-1}(G)$ is equal to its Kirchhoff index,
denoted by $Kf(G)$, which found applications in electric circuit,
probabilistic theory and chemistry \cite{GuMo,Pal}. Some properties
of $s_{\alpha}$ for $\alpha\ne 0, 1$, including further properties
of $s_{2}$ and $s_{\frac{1}{2}}$ have been established recently in
\cite{Zhou}. Now we give further properties of $s_{\alpha}$, that
is, bounds related to the degree sequences of the graphs. As a
by-product, a lower bound for the Kirchhoff index is given.

Note that lots of spectral indices were proposed in \cite{VT}
recently, and since the Laplacian eigenvalues are all nonnegative,
for $\alpha\ne 0$, $s_{\alpha}$ is equal to the spectral index
$SpSum^{\alpha}(L)$ with $L$ being the Laplacian matrix of the
graph.

The Estrada index of a graph $G$ with eigenvalues $\lambda_1,
\lambda_2, \dots, \lambda_n$ is defined as
$EE(G)=\sum_{i=1}^ne^{\lambda_i}$.  It is a very useful descriptors
in a large variety of problems, including those in biochemistry and
in complex networks [9--11],
for recent results see [12--14]. 
The Laplacian Estrada index of a graph $G$ with $n$ vertices is
defined as \cite{GFA}
\[
LEE(G)=\sum_{i=1}^n e^{\mu_i}.
\]
We also give bounds for the Laplacian Estrada index related to the
degree sequences of the graphs.

\vspace{4mm}

\begin{center}
{\bf 2. PRELIMINARIES}
\end{center}

For two non-increasing sequences $x=(x_{1}, x_{2}, \dots, x_{n})$
and $y=(y_{1}, y_{2}, \dots, y_{n})$,  $x$ is majorized by $y$,
denoted by $x\preceq y$,  if
\begin{eqnarray*}
&&\sum_{i=1}^{j}x_{i}\leq\sum_{i=1}^{j}y_{i} \mbox{ for } j=1, 2,
\dots,n-1, \mbox{ and}\\
&&\sum_{i=1}^{n}x_{i}=\sum_{i=1}^{n}y_{i}.
\end{eqnarray*}
For a real-valued function $f$ defined on a set in $\mathbb{R}^n$,
if $f(x)<f(y)$ whenever $x\preceq y$ but $x\ne y$, then $f$ is said
to be strictly Schur-convex \cite{MaOl}.


\begin{Lemma} \label{schur} Let $\alpha$ be a real number with  $\alpha\ne
0,1$.

$(i)$ For $x_i\ge 0$, $i=1,2, \dots, h$,  $f(x)=\sum_{i=1}^h
x_i^{\alpha}$ is strictly Schur-convex if  $\alpha>1$, and
$f(x)=-\sum_{i=1}^h x_i^{\alpha}$ is strictly Schur-convex if
$0<\alpha<1$.

$(ii)$ For $x_i>0$, $i=1,2, \dots, h$, $f(x)=\sum_{i=1}^h
x_i^{\alpha}$ is strictly Schur-convex if $\alpha<0$.
\end{Lemma}

\begin{Proof} From \cite[p. 64, C.1.a]{MaOl} we know that if the real-valued function
$g$ defined on an interval in $\mathbb{R}$ is a strictly convex then
$\sum_{i=1}^ hg(x_i)$ is strictly Schur-convex.

If $x_i\ge 0$, then $x_i^{\alpha}$ is strictly convex if $\alpha>1$
and $-x_i^{\alpha}$ is strictly convex if $0<\alpha<1$, and thus (i)
follows.

If $x_i>0$ and $\alpha<0$, then $x_i^{\alpha}$ is strictly convex,
and thus (ii) follows.
\end{Proof}

Let $K_n$ and $S_n$ be respectively the complete graphs and the star
with $n$ vertices.  Let $K_n-e$ be the graph with one edge deleted
from $K_n$.

Recall the the degree sequence of a graph $G$ is a list of the
degrees of the vertices in  non-increasing order, denoted by $(d_1,
d_2, \dots, d_n)$, where $n$ is the number of vertices of $G$.  Then
$d_1$ is the maximum vertex degree of $G$.

\vspace{4mm}

\begin{center}
{\bf 3. BOUNDS FOR $s_{\alpha}$ RELATED TO DEGREE SEQUENCES}
\end{center}

 We need the following lemmas.

\begin{Lemma}\label{maj} \textnormal{\cite{Gro}} Let $G$ be a connected graph with
$n\ge 2$ vertices.  Then $(d_1+1, d_2, \dots, d_{n-1}, d_n-1)\preceq
(\mu_1, \mu_2, \ldots, \mu_n)$.
\end{Lemma}

\begin{Lemma} \label{known} \textnormal{\cite{Zhou}}
\label{start2} Let $G$ be a connected graph with $n\ge 2$ vertices.
Then $\mu_2=\cdots=\mu_{n-1}$ and $\mu_1=1+d_1$ if and only if
$G=K_n$ or $G=S_n$.
\end{Lemma}

Now we provide bounds for $s_{\alpha}$ using degree sequences.

\begin{Proposition} \label{Prop1} Let $G$ be a connected graph with $n\ge 2$
vertices. Then
\begin{eqnarray*}
&&s_{\alpha}(G)\ge (d_1+1)^{\alpha}+\sum_{i=2}^{n-1}
d_i^{\alpha}+(d_n-1)^{\alpha} \ if \ \alpha>1\\
&& s_{\alpha}(G)\le (d_1+1)^{\alpha}+\sum_{i=2}^{n-1}
d_i^{\alpha}+(d_n-1)^{\alpha} \ if \ 0<\alpha<1
\end{eqnarray*}
with either equality if and only if $G=S_n$.
\end{Proposition}

\begin{Proof}
If $\alpha>1$, then by Lemma~\ref{schur} (i), $f(x)=\sum_{i=1}^n
x_i^{\alpha}$ is strictly Schur-convex, which, together with
Lemma~\ref{maj}, implies that
\[
s_{\alpha}(G)=\sum_{i=1}^n \mu_i^{\alpha}\ge
(d_1+1)^{\alpha}+\sum_{i=2}^{n-1} d_i^{\alpha}+(d_n-1)^{\alpha}
\]
with equality if and only if $(\mu_1, \mu_2, \ldots, \mu_n)=(d_1+1,
d_2, \dots, d_{n-1}, d_n-1)$.

If $0<\alpha<1$, then by Lemma~\ref{schur} (i), $f(x)=-\sum_{i=1}^h
x_i^{\alpha}$ is strictly Schur-convex, which, together with
Lemma~\ref{maj}, implies that
\[
-s_{\alpha}(G)=-\sum_{i=1}^n \mu_i^{\alpha}\ge
-\left[(d_1+1)^{\alpha}+\sum_{i=2}^{n-1}
d_i^{\alpha}+(d_n-1)^{\alpha}\right],
\]
i.e.,
\[
s_{\alpha}(G)=\sum_{i=1}^n \mu_i^{\alpha}\le
(d_1+1)^{\alpha}+\sum_{i=2}^{n-1} d_i^{\alpha}+(d_n-1)^{\alpha}
\]
with equality if and only if $(\mu_1, \mu_2, \ldots, \mu_n)=(d_1+1,
d_2, \dots, d_{n-1}, d_n-1)$.

By Lemma~\ref{known}, we have $(\mu_1, \mu_2, \ldots, \mu_n)=(d_1+1,
d_2, \dots, d_{n-1}, d_n-1)$  if and only if  $G=S_{n}$.
\end{Proof}

We note that the result  for  $\alpha=\frac{1}{2}$ has been given in
\cite{LL}.

\begin{Proposition} \label{Prop2} Let $G$ be a connected graph with $n\ge 3$
vertices. If $\alpha<0$, then
\begin{eqnarray*}
s_{\alpha}(G)\ge (d_1+1)^{\alpha}+\sum_{i=2}^{n-2}
d_i^{\alpha}+\left(d_{n-1}+d_n-1\right)^{\alpha}
\end{eqnarray*}
with equality if and only if $G=S_n$ or $G=K_3$.
\end{Proposition}

\begin{Proof}
By Lemma~\ref{schur} (ii), $f(x)=\sum_{i=1}^{n-1}x_i^{\alpha}$ is
strictly Schur-convex for $x_i>0$, $i=1,2,\dots, n-1$.  By
Lemma~\ref{maj}, $(d_1+1, d_2, \dots, d_{n-2}, d_{n-1}+d_n-1)\preceq
(\mu_1, \mu_2, \ldots, \mu_{n-1})$. Thus
\[
s_{\alpha}(G)=\sum_{i=1}^{n-1} \mu_i^{\alpha}\ge
(d_1+1)^{\alpha}+\sum_{i=2}^{n-2}
d_i^{\alpha}+\left(d_{n-1}+d_n-1\right)^{\alpha}
\]
with equality if and only if $(\mu_1, \mu_2, \ldots,
\mu_{n-1})=(d_1+1, d_2, \dots, d_{n-2}, d_{n-1}+d_n-1)$, which, by
Lemma~\ref{known}, is equivalent to $G=S_{n}$ or $G=K_3$.
\end{Proof}

Let $G$ be a connected graph with $n\ge 3$ vertices. Then by
Proposition \ref{Prop2},
\begin{eqnarray*}
Kf(G)\ge
n\left(\frac{1}{d_1+1}+\sum_{i=2}^{n-2}\frac{1}{d_i}+\frac{1}{d_{n-1}+d_n-1}\right)
\end{eqnarray*}
with equality if and only if $G=S_n$ or $G=K_3$. Note that we have
already shown in \cite{ZhTr} that
\[
Kf(G)\ge -1+(n-1)\sum_{i=1}^n\frac{1}{d_i}.
\]
These two lower bounds are incomparable  as for $K_n$ with $n\ge 4$
the latter is better but for $K_n-e$ with $n\ge 7$ the former is
better.

\vspace{3mm}

\noindent {\bf Remark 1.} For the degree sequence $(d_{1}, d_{2},
\dots, d_{n})$ of a graph, its conjugate sequence is
$(d_{1}^*,d_{2}^*, \dots, d_{n}^*)$, where $d_i^*$ is equal to the
cardinality of the set $\{j:d_j\ge i\}$. Note that $(d_1, d_2,
\dots, d_n)\preceq (d_1^*, d_2^*, \dots, d_n^*)$ \cite{Mer,MM}. It
was conjectured in \cite{MM} that
\[
(\mu_1, \mu_2, \dots, \mu_n)\preceq (d_1^*, d_2^*, \dots, d_n^*).
\]
Though  still open, it has been proven to be true for a class of
graphs including trees \cite{Ste}. Let $G$ be a tree with $n\ge 2$
vertices. Then $d_1^*=n$, $d_{d_1+1}^*=0$,  and by similar arguments
as in the proof of Proposition \ref{Prop1}, we have
\begin{eqnarray*}
&&s_{\alpha}(G)\le \sum_{i=1}^{d_1}
\left(d_i^*\right)^{\alpha} \mbox{ if } \alpha>1 \mbox{ or } \alpha<0\\
&& s_{\alpha}(G)\ge \sum_{i=1}^{d_1} \left(d_i^*\right)^{\alpha}
\mbox{ if } 0<\alpha<1
\end{eqnarray*}
with either equality if and only if $(\mu_1, \mu_2, \ldots,
\mu_n)=(d_1^*, d_2^*, \dots, d_n^*)$, which,  is equivalent to
$G=S_n$ since if $G\ne S_n$, then $d_{n-1}^*=0$ but $\mu_{n-1}>0$.

To end this section, we mention a result of Rodriguez and Petingi
concerning the Laplacian spectral moments in \cite{RoP}:

\begin{Proposition}\label{RP}
For a graph $G$ with $n$ vertices and any positive integer $k$, we
have
\[
s_{k}(G)\ge \sum_{i=1}^nd_i(1+d_i)^{k-1}
\]
and for $k\ge 3$, equality occurs if and only if $G$ is a
vertex--disjoint union of complete subgraphs.
\end{Proposition}

\vspace{4mm}

\begin{center}
{\bf 4. BOUNDS FOR LAPLACIAN ESTRADA INDEX RELATED TO DEGREE
SEQUENCES}
\end{center}

Let $G$ be a graph with $n$ vertices.
Obviously,
\[
LEE(G)=\sum_{k\ge 0}\frac{t_{k}(G)}{k!}=n+\sum_{k\ge
1}\frac{s_{k}(G)}{k!}.
\]
Thus, properties of the Laplacian moments in previous section may be
converted into properties of the Laplacian Estrada index.

\begin{Proposition} \label{lee1} Let $G$ be a connected graph with $n\ge 2$
vertices. Then
\[ LEE(G)\ge
e^{d_1+1}+\sum_{i=2}^{n-1}e^{d_i}+e^{d_n-1}
\]
with equality if and only if $G=S_n$.
\end{Proposition}

\begin{Proof} Note that $t_0(G)=n$, $t_1(G)=\sum_{i=1}^nd_i$,  and $t_k(G)=s_k(G)$ for $k\ge 1$. By Proposition
\ref{Prop1},
\[
t_k(G)\ge (d_1+1)^{k}+\sum_{i=2}^{n-1} d_i^{k}+(d_n-1)^{k}
\]
for $k=0, 1, \dots$, with equality for $k=0, 1$, and  if $k\ge 2$
then  equality occurs if and only if $G=S_n$.  Thus
\begin{eqnarray*}
LEE(G)&=&\sum_{k\ge 0}\frac{t_{k}(G)}{k!}\\
&\ge & \sum_{k\ge
0}\frac{(d_1+1)^{k}+\sum_{i=2}^{n-1} d_i^{k}+(d_n-1)^{k}}{k!}\\
&=&e^{d_1+1}+\sum_{i=2}^{n-1}e^{d_i}+e^{d_n-1}
\end{eqnarray*}
with equality if and only if $G=S_n$.
\end{Proof}

Similarly, if $G$ be a tree with $n\ge 2$ vertices, Then 
by similar arguments as in the proof of Proposition \ref{lee1},
we have
\[
LEE(G)\le \sum_{i=1}^ne^{d_i^*}=n-d_1+\sum_{i=1}^{d_1}e^{d_i^*}
\]
with equality if and only if  $G=S_n$.

\begin{Proposition} Let $G$ be a graph with $n\ge 2$
vertices. Then
\[ LEE(G)\ge
n+\sum_{i=1}^n\frac{d_i}{1+d_i}\left(e^{1+d_i}-1\right)
\]
with equality if and only if $G$ is a vertex--disjoint union of
complete subgraphs.
\end{Proposition}

\begin{Proof} By Proposition \ref{RP},
\[
t_k(G)\ge \sum_{i=1}^nd_i(1+d_i)^{k-1}
\]
for $k=1, 2\dots$, and for $k\ge 3$ equality occurs if and only if
$G$ is a disjoint union of cliques. The inequality above is an
equality for $k=1,2$.  Thus
\begin{eqnarray*}
LEE(G)&=&\sum_{k\ge 0}\frac{t_{k}(G)}{k!}\\
&\ge & n+\sum_{k\ge
1}\frac{\sum_{i=1}^nd_i(1+d_i)^{k-1}}{k!}\\
&=& n+\sum_{i=1}^n\frac{d_i}{1+d_i}\sum_{k\ge 1}\frac{(1+d_i)^{k}}{k!}\\
&=& n+\sum_{i=1}^n\frac{d_i}{1+d_i} \left(e^{1+d_i}-1\right)
\end{eqnarray*}
with equality if and only if $G$ is a vertex--disjoint union of
complete subgraphs.
\end{Proof}

\vspace{3mm}

\noindent {\bf Remark 2.} We note that lower bounds on the Laplacian
spectral moments in \cite{Zhou} may also be converted to the bounds
of Laplacian Estrada index.

(a) Let $G$ be a connected graph with $n\ge 3$ vertices, $m$ edges.
 Then
\[
LEE(G)\ge 1+e^{1+d_1}+(n-2) e^{\frac{2m-1-d_1}{n-2}}
\]
\[
LEE(G)\ge 1+e^{1+d_1}+(n-2)
e^{\left(\frac{tn}{1+d_1}\right)^{\frac{1}{n-2}}}
\]
with either equality if and only if $G=K_n$ or $G=S_n$, where $t$ is
the number of spanning trees in $G$.

(b) Let $G$ be a graph with $n\ge 2$ vertices and $m$ edges. Let
$\overline{G}$ be the complement of the graph $G$. By the
arithmetic--geometric inequality, we have
$LEE(G)=1+\sum\limits_{i=1}^{n-1}e^{\mu_i}\ge 1+(n-1)
e^{\frac{2m}{n-1}}$ with equality if and only if
$\mu_1=\mu_2=\cdots=\mu_{n-1}$, i.e., $G=K_n$ or $G=\overline{K_n}$
\cite{Zhou}. Let $\overline{m}$ be the number of edges of
$\overline{G}$\,. Thus
\begin{eqnarray*}
LEE(G)+LEE(\overline{G})&\ge & 2+(n-1)\left(
e^{\frac{2m}{n-1}}+e^{\frac{2\overline{m}}{n-1}}\right)\\
&\ge & 2+2(n-1)e^{\frac{2m+2\overline{m}}{2(n-1)}}\\
&=& 2+2(n-1)e^{\frac{n}{2}},
\end{eqnarray*}
and then $LEE(G)+LEE(\overline{G})> 2+2(n-1)e^{\frac{n}{2}}$.

(c)   Let $G$ be a connected bipartite graph with $n\ge 3$ vertices
and $m$ edges.  Recall that the first Zagreb index of a graph $G$,
denoted by $M_1(G)$, is defined as the sum of the squares of the
degrees of the graph [22--24]. 
Then
\[
LEE(G)\ge 1+e^{2\sqrt{\frac{M_1(G)}{n}}}+(n-2)
e^{\frac{2m-2\sqrt{\frac{M_1(G)}{n}}}{n-2}}
\]
\[
LEE(G)\ge 1+e^{2\sqrt{\frac{M_1(G)}{n}}}+(n-2)
e^{\left(\frac{tn\sqrt{n}}{2\sqrt{M_1(G)}}\right)^{\frac{1}{n-2}}}
\]
with either equality if and only if  $n$ is even and
$G=K_{\frac{n}{2},\frac{n}{2}}$, where $t$ is the number of spanning
trees in $G$.

\vspace{6mm}

\baselineskip=0.25in

\noindent {\it Acknowledgement.\/} This work was supported by the
National Natural Science Foundation  (no. 10671076) and the
Guangdong Provincial Natural Science Foundation  (no.
8151063101000026) of China.

\end{document}